\documentclass[12pt]{article}
\usepackage{amsmath,amsfonts,amssymb}
\usepackage{color}
\def\C{\mathbb{C}}
\def\N{\mathbb{N}}
\def\R{\mathbb{R}}
\def\a{\mathbf{a}}
\def\b{\mathbf{b}}
\def\L{\mathcal{L}}
\definecolor{greenish}{rgb}{0.1,0.7,0}
\newtheorem{theorem}{\hspace*{\parindent}Theorem}
\newtheorem{lemma}{\hspace*{\parindent}Lemma}
\newtheorem{corollary}{\hspace*{\parindent}Corollary}
\title{Applications of the Stieltjes and Laplace transform representations of the hypergeometric functions}
\author{D.B.\:Karp$^{\rm a,b}$\footnote{Corresponding author. E-mail: D.B.\:Karp -- \emph{dimkrp@gmail.com}, E.G.Prilepkina --  \emph{pril-elena@yandex.ru}}~~and E.G.Prilepkina$^{\rm a,b}$
\\[10pt]
\small{\textit{$\phantom{1}^a$Far Eastern Federal University, Vladivostok, Russia}}\\\small{\textit{$\phantom{1}^b$Institute of Applied Mathematics, FEBRAS}}}
\date{}
\begin{document}
\maketitle

\begin{abstract}
In our previous work we found sufficient conditions to be imposed on the parameters of the generalized hypergeometric function in order that it be
completely monotonic or of Stieltjes class.  In this paper we collect a number of consequences of these properties.  In particular, we find new integral representations of the generalized hypergeometric functions, evaluate a number of integrals of their products, compute the jump and the average value of the the generalized hypergeometric function over the branch cut $[1,\infty)$, establish new inequalities for this function in the half plane $\Re{z}<1$.  Furthermore, we discuss integral representations of absolutely monotonic functions and present a curious formula for a finite sum of products of gamma ratios as an integral of Meijer's $G$ function.
\end{abstract}

\bigskip


\bigskip


\bigskip

\section{Introduction and preliminaries}

Throughout the history of the hypergeometric functions integral representations played an important role in their study. The celebrated Euler's integral for the Gauss functions ${}_2F_{1}$ \cite[Theorem~2.2.1]{AAR} was probably the first among them.  More recently, Kiryakova  \cite[Chapter~4]{KiryakovaBook} observed that the generalized hypergeometric functions possess similar representations with densities expressed via Meijer's $G$ function.  These representations have been extended and massively applied in our papers \cite{KarpJMS,KLMeijer1,KPJMAA,KSJAT}.  For the generalized hypergeometric functions of the Gauss and Kummer type the most important integral representations are those given by the (generalized) Stieltjes and Laplace transforms, respectively.  If the density in the corresponding representation is nonnegative,  then the resulting hypergeometric function is completely monotonic or belongs to the (generalized) Stieltjes class.

The purpose of this note is to share some observations about completely monotonic and Stieltjes functions and illustrate them by the hypergeometric examples.  As a by-product of these observations we evaluate a number of integrals involving the generalized hypergeometric functions which are neither contained in the most comprehensive collection \cite{PBM3} nor  are solvable by \textit{Mathematica}. We further find the jump and the average value of the generalized hypergeometric function over the branch cut and present some inequalities for this function resulting from the Stieltjes representation.  Some facts presented in this note complement and/or illustrate the results of the recent work \cite{KoumPeders} by Koumandos and Pedersen.

We now turn to the details.  Let us fix some notation and terminology first.  The standard symbols $\N$, $\R$ and $\C$ will be used to denote the natural, real and complex numbers, respectively; $\N_0=\N\cup\{0\}$, $\R_{+}=[0,\infty)$.  A nonnegative function $f$ defined on $(0,\infty)$ is called completely monotonic if it has derivatives of all orders and $(-1)^nf^{(n)}(x)\ge0$ for $n\ge0$ and $x>0$. A function $\phi:[0,\infty)\mapsto\R$ is called absolutely monotonic if it is infinitely differentiable on $[0,\infty)$ and $\phi^{(k)}(x)\ge0$ for all $k\ge0$ and all $x\ge0$. According to \cite[Theorem~3a]{Widder} an absolutely monotonic function $\phi$ on $[0,\infty)$ has an extension to an entire function with the power series expansion $\phi(z)=\sum_{n\ge0}\phi_nz^n$, where $\phi_n\ge0$ for all $n\ge0$.  As usual $\Gamma(z)$ stands for Euler's gamma function and $(a)_n=\Gamma(a+n)/\Gamma(a)$ denotes the rising factorial (or Pochhammer's symbol).  Throughout the paper we will use the shorthand notation for the products and sums. For $\a\in\C^p$ define
\begin{equation*}
\begin{split}
&\Gamma(\a)=\Gamma(a_1)\Gamma(a_2)\cdots\Gamma(a_p),~~~(\a)_n=(a_1)_n(a_2)_n\cdots(a_p)_n,
\\
&\a+\mu=(a_1+\mu,a_2+\mu,\dots,a_p+\mu),~~~~\a>0~\Leftrightarrow~a_k>0~\text{for}~k=1,\ldots,p;
\end{split}
\end{equation*}
$(\a)_1$ we be further simplified to $(\a)=\prod_{k=1}^{p}a_k$. We will write $\a_{[k]}$ for the vector $\a$ with $k$-th element removed, i.e. $\a_{[k]}=(a_1,\ldots,a_{k-1},a_{k+1},\ldots,a_p)$. Further, we will follow the standard definition of the generalized hypergeometric function ${_{p}F_q}$ as the sum of the series
\begin{equation}\label{eq:pFqdefined}
{_{p}F_q}\left(\left.\!\!\begin{array}{c}\a \\ \b\end{array}\right|z\!\right)={_{p}F_q}\left(\a;\b;z\right)
=\sum\limits_{n=0}^{\infty}\frac{(a_1)_n(a_2)_n\cdots(a_{p})_n}{(b_1)_n(b_2)_n\cdots(b_q)_nn!}z^n
\end{equation}
if $p\le{q}$, $z\in\C$. If $p=q+1$ the above series has unit radius of convergence and ${_{p}F_q}(z)$ is defined as analytic continuation of its sum to $\C\setminus[1,\infty)$.  Here $\a=(a_1,\ldots,a_p)$, $\b=(b_1,\ldots,b_q)$ are (generally complex) parameter vectors, such that $-b_j\notin\N_0$, $j=1,\ldots,q$.

According to \cite[Corollary~1]{KarpJMS}
\begin{equation}\label{eq:pFpLaplace}
{_{p}F_p}\!\left(\left.\!\!\begin{array}{c}\a\\
\b\end{array}\right|-z\!\right)
=\frac{\Gamma(\b)}{\Gamma(\a)}\int\nolimits_{0}^{1}e^{-zt}G^{p,0}_{p,p}\left(t\left|\begin{array}{l}\!\!\b\!\!\\\!\!\a\!\!\end{array}\right.\right)
\frac{dt}{t}
\end{equation}
and
\begin{equation}\label{eq:p+1FpLaplace}
{_{p+1}F_p}\!\left(\left.\!\!\begin{array}{c}\a\\
\b\end{array}\right|-z\!\right)
=\frac{\Gamma(\b)}{\Gamma(\a)}\int\nolimits_{0}^{\infty}e^{-zt}G^{p+1,0}_{p,p+1}\left(t\left|\begin{array}{l}\!\!\b\!\!\\\!\!\a\!\!\end{array}\right.\right)
\frac{dt}{t}
\end{equation}
under convergence conditions of the integrals on the right hand side, as elucidated in \cite[Theorem~1]{KarpJMS}.  Here $G^{m,n}_{p,q}$ denotes Meijer's $G$
function whose definition and properties can be found in  \cite[Section~16.17]{NIST}, \cite[Section~12.3]{BealsWong}, \cite[Section~8.2]{PBM3},  \cite[\S2]{KarpJMS} and \cite[Section~2]{KPJMAA}.  Formulas (\ref{eq:pFpLaplace}) and (\ref{eq:p+1FpLaplace}) combined with Bernstein's theorem imply that the hypergeometric functions on the left hand sides are completely monotonic if the $G$ functions in the integrands are nonnegative.  According to \cite[Theorem~2]{KarpJMS} the function $G^{p,0}_{p,p}$ from (\ref{eq:pFpLaplace}) is nonnegative if the \emph{M\"{u}ntz polynomial}
\begin{equation}\label{eq:v-defined}
v_{\a,\b}(t)=\sum\limits_{k=1}^{p}(t^{a_k}-t^{b_k})\ge0~\text{for}~t\in[0,1].
\end{equation}
Note that the vector $\a$ in (\ref{eq:p+1FpLaplace}) contains $p+1$ element, while $v_{\a,\b}(t)$ is defined for two real vectors $\a$, $\b$ of equal size. Therefore, $v_{\a,\b}(t)$ in (\ref{eq:v-defined}) should be replaced by $v_{\a_{[p+1]},\b}(t)$ for $G^{p+1,0}_{p,p+1}$ to be nonnegative, while the remaining parameter $a_{p+1}$ may take arbitrary positive values.  See \cite[Proof of Th.~8]{KLMeijer1} for details.
Inequality (\ref{eq:v-defined}) is implied by the stronger condition $\b\prec^W\a$ known as the weak supermajorization \cite[section~2]{KPCMFT}  and given by  \cite[Definition~A.2]{MOA}
\begin{equation}\label{eq:amajorb}
\begin{split}
& 0<a_1\leq{a_2}\leq\cdots\leq{a_p},~~
0<b_1\leq{b_2}\leq\cdots\leq{b_p},
\\
&\sum\limits_{i=1}^{k}a_i\leq\sum\limits_{i=1}^{k}b_i~~\text{for}~~k=1,2\ldots,p.
\end{split}
\end{equation}
Further sufficient conditions for the validity of (\ref{eq:v-defined}) in terms of $\a$, $\b$ can be found in our recent paper \cite[section~2]{KPCMFT}.

This paper is organized as follows.  Section~2 is concerned with the Laplace transform representations and their corollaries, such as integral representations of absolutely monotonic functions and new evaluations of convolution integrals containing hypergeometric functions. Sections~3 is devoted to the consequences of the (generalized) Stieltjes transform representation for the Gauss type generalized hypergeometric functions. These are: new integrals for ${}_{p+1}F_{p}(z)$ in the left half plane and the sector $\arg(-z)\in(-\pi/3,\pi/3)$, a formula for the jump of this functions over the branch cut $(-\infty,-1]$ and for its mean on the banks of the cut, inequalities for this function in the half-plane $\Re(z)<1$ implied by its univalence. Finally, in the ultimate Section~4 we present a curious formula for the finite sum of the differences of ratios of the gamma products as an integral of $G$ function over the interval $(0,1)$.

\section{Laplace transform representations and their consequences}

According to \cite[Theorem~4]{KarpJMS}, the function
$$
x\mapsto x^{-a_{p+1}}{_{p+1}F_p}\!\left(\left.\!\!\begin{array}{c}\a\\
\b\end{array}\right|-1/x\!\right)
$$
is completely monotonic if $\a>0$, $v_{\a_{[p+1]},\b}(t)\ge0$ for $t\in[0,1]$ and $a_{p+1}>0$. Using our recent result from \cite{KLMeijer1} we can prove a similar statement for ${}_pF_p$. Recall that $\a_{[k_1,k_2,\ldots,k_r]}$ denotes the vector $\a$ with the elements $a_{k_1},a_{k_2},\ldots,a_{k_r}$ removed.
\begin{theorem}\label{th:pFpnewcase}
Suppose $p\ge2$, $\a,\b\in\R^{p}$ are ascending positive vectors and set $\a_*=(\a_{[1,p]},3/2)$ and $\b_*=\b_{[1]}$. Assume that
$0<a_1\le1$, $b_1\ge{a_1+1}$ and $v_{\a_*,\b_*}(t)\ge0$ on $[0,1]$.  Then the function
$$
x\mapsto x^{-a_p}{_{p}F_{p}}(\a;\b;-1/x)
$$
is completely monotonic.
\end{theorem}
\textbf{Proof.}  Termwise integration confirms that
$$
x^{-a_p}{_{p}F_{p}}(\a;\b;-1/x)=\frac{1}{\Gamma(a_p)}\!\int\limits_{0}^{\infty}\!e^{-xt}t^{a_p-1}{}_{p-1}F_p(\a_{[p]};\b;-t)dt.
$$
Under conditions of the theorem the function  ${}_{p-1}F_p(\a_{[p]};\b;-t)$ is positive on $[0,\infty)$ by \cite[Theorem~7]{KLMeijer1} and the claim follows by Bernstein's theorem. $\hfill\square$

It is immediate from the definition of a completely monotonic function that a linear combination of such functions with nonnegative coefficients is completely monotonic.  As pointwise limit of a sequence of completely monotonic functions is also completely monotonic \cite[p.151]{Widder}, we conclude that $\phi(1/x)$ is completely monotonic if $\phi(x)$ is absolutely monotonic (see also \cite[Theorem~3]{MS}).  We summarize these facts in the next simple observation.
\begin{theorem}\label{th:twoCMfunctions}
Suppose $\phi:[0,\infty)\mapsto\R$ is absolutely monotonic. Then  $\phi(1/x)$ is completely monotonic and its representing measure is given by
$\nu(dt)=[\phi_0\delta_0+\nu(t)]dt$, where $\nu(t)=\sum_{k\ge0}\phi_{k+1}t^k/k!$.  If, furthermore, $\{k!\phi_{k}\}_{k\ge0}$ forms a Stieltjes moment sequence, then
$\phi(-x)$ is also completely monotonic and its representing measure $\mu(dt)$ is related to $\nu(t)$ by
$$
\nu(t)=\int\limits_{[0,\infty)}\!\!\left(\frac{u}{t}\right)^{1/2}\!I_1(2\sqrt{ut})\mu(du),
$$
where $I_1$ is the modified Bessel function of the first kind.
\end{theorem}
\textbf{Proof.} Indeed, according to \cite[Theorem~2b]{Widder} the function $x\mapsto\phi(1/x)$ is completely monotonic. Furthermore,
\begin{multline*}
\int\limits_{[0,\infty)}e^{-xt}[\phi_0\delta_0+\nu(t)]dt=
\phi_0+\sum_{k\ge0}\frac{\phi_{k+1}}{k!}\int\limits_{[0,\infty)}e^{-xt}t^kdt
\\
=\phi_0+\sum_{k\ge0}\frac{\phi_{k+1}}{k!}\frac{\Gamma(k+1)}{x^{k+1}}
=\sum_{k\ge1}\frac{\phi_{k}}{x^{k}}+\phi_0=\phi(1/x),
\end{multline*}
which proves the first claim. On the other hand, by the Stieltjes moment property $k!\phi_k=\int\limits_{[0,\infty)}t^k\mu(dt)$ for some nonnegative measure $\mu$, so that
$$
\phi(-x)=\sum\limits_{k=0}^{\infty}\phi_k(-x)^k
=\sum\limits_{k=0}^{\infty}\frac{(-x)^k}{k!}\int\limits_{[0,\infty)}t^k\mu(dt)
=\int\limits_{[0,\infty)}e^{-xt}\mu(dt).
$$
Then,
\begin{multline*}
\nu(t)=\sum_{k\ge0}\phi_{k+1}t^k/k!=\sum_{k\ge0}\frac{t^k}{k!(k+1)!}\int\limits_{[0,\infty)}u^{k+1}\mu(du)
\\
=\int\limits_{[0,\infty)}u\sum_{k\ge0}\frac{(ut)^k}{(2)_kk!}\mu(du)
=\int\limits_{[0,\infty)}u{}_0F_1(-;2;ut)\mu(du)
=\int\limits_{[0,\infty)}\!\!\left(\frac{u}{t}\right)^{1/2}\!I_1(2\sqrt{ut})\mu(du),
\end{multline*}
where we have used the series expansion for the modified Bessel function \cite[formula (10.25.2)]{NIST}.
$\hfill\square$

\smallskip
\textbf{Remark.}  Theorem~\ref{th:twoCMfunctions} shows that the  absolutely monotonic functions on $[0,\infty)$ are precisely those that are representable in the form
$$
\phi(y)=C+\int_{0}^{\infty}e^{-t/y}\nu(t)dt,
$$
where $C\ge0$ and $\nu(t)$ is an entire function of minimal exponential type (i.e. of order $<1$ or of order $=1$ and minimal type) with nonnegative Taylor coefficients $\nu_k$. The minimal exponential type simply means that $\limsup_{k\to\infty}[\nu_kk!]^{1/k}=0$ \cite[section~1.3]{Levin}.  Other integral representations can be obtained using \cite[Theorem~2b]{Widder} which asserts that $\phi(f(x))$ is completely monotonic once $\phi$ is absolutely monotonic and $f(x)$ is  completely monotonic.  The class of the representing measures will depend on $f(x)$. Taking, for instance, $f(x)=e^{-x}$ yields $\phi(y)=\int_{[0,\infty)}y^t\mu(dt)$.  This is, of course, just another way of writing the Taylor series expansion, so that $\mu(dt)$ must be supported on nonnegative integers with mass concentrated at $n$ decreasing fast enough.  Note, that the problem of finding an integral representation for absolutely monotonic functions has been considered in \cite{Sitnik}.

\smallskip
\textbf{Remark.} The density $\nu(t)$ of the representing measure for $\phi(1/x)$ in Theorem~\ref{th:twoCMfunctions} is also absolutely monotonic. Laplace transforms of absolutely monotonic functions have been recently studied by Koumandos and Pedersen in \cite[Theorem~1.1]{KoumPeders}.  The proof of their theorem indicates  that the authors essentially consider the functions of the form $(1/x)\phi(1/x)$. This shows that Theorem~\ref{th:twoCMfunctions} complements \cite[Theorem~1.1]{KoumPeders}.

\smallskip
According to Theorem~\ref{th:twoCMfunctions} the function $x\mapsto{_{p}F_q}(\a;\b;1/x)$ is completely monotonic for $p\le{q}$ and any positive parameter vectors $\a$, $\b$. We write the representing measure explicitly in the next corollary.
\begin{corollary}\label{cr:pFq1-over-x}
Suppose $p\leq{q+1}$,  $\a$ and $\b$ are any complex vectors, $\b$ does not contain non-positive integers. Then the following identity holds:
\begin{equation}\label{eq:pFq1-over-x}
f(x)={_{p}F_q}\left(\a;\b;1/x\right)=\int\limits_{[0,\infty)}e^{-xt}\left[\frac{(\a)}{(\b)}{_{p}F_{q+1}}\left(\left.\!\!\begin{array}{c}\a+1\\\b+1,2\end{array}\right|t\!\right)+\delta_0\right]dt,
\end{equation}
where $\delta_{0}$ is the Dirac measure with mass $1$ concentrated at zero.  Here $x>0$ if $p\le{q}$ and $x>1$ if $p=q+1$.
\end{corollary}
\textbf{Proof.} Straightforward termwise integration yields the claim. $\hfill\square$

Curiously enough, formula (\ref{eq:pFq1-over-x}) is different from the related integral evaluations \cite[Exercise~11, p.115]{AAR} and \cite[(2.22.3.1)]{PBM3}.

Further, by application of \cite[Theorem~1.1]{KoumPeders} to the function $f$ defined in (\ref{eq:pFq1-over-x}) the derivatives $(-1)^k(x^kf(x))^{(k)}$  are completely monotonic for each $k\in\N_0$.  Moreover, if $R_n(y)$ is $n$-th Taylor remainder for ${_{p}F_q}(\a;\b;y)$ then $R_n(1/x)$ is completely monotonic of order $n$ ($g:(0,\infty)\mapsto\R$ is completely monotonic of order $\alpha>0$ if $x\to{x^{\alpha}g(x)}$ is completely monotonic). These properties are relatively straightforward in the case of the above $f$, but they may serve as a good illustration for
\cite[Theorem~1.1]{KoumPeders}.  A slightly different hypergeometric illustration is given in \cite[Remark~4.3]{KoumPeders}.

\textbf{Remark.} From what we have said both ${_{p}F_p}(\a;\b;-x)$  and  ${_{p}F_p}(\a;\b;1/x)$ are completely monotonic if we assume that $v_{\a,\b}(t)=\sum_{j=1}^{p}(t^{a_j}-t^{b_j})\ge0$ for $t\in[0,1]$.  If, in addition, conditions of Theorem~\ref{th:pFpnewcase} are satisfied, then all three functions
$$
x\mapsto{_{p}F_p}(\a;\b;-x),~~~~x\mapsto{_{p}F_{p}}(\a;\b;1/x),~~~~x\mapsto x^{-a_p}{_{p}F_{p}}(\a;\b;-1/x)
$$
are completely monotonic.

\smallskip

Next, we want to combine (\ref{eq:pFq1-over-x}) with the product formulas for hypergeometric functions (conveniently collected, for instance, in \cite{Grinshpan}) and employ Laplace convolution theorem to get new integral evaluations for hypergeometric functions.  First we handle the general case. Suppose we have an identity of the form $\phi_1(y)\phi_2(y)=\phi_0(y)$, where all three functions are absolutely monotonic and the representing measures in Theorem~\ref{th:twoCMfunctions} are $\nu_1(t)dt+\delta_0$, $\nu_2(t)dt+\delta_0$ and $\nu_0(t)dt+\delta_0$, respectively (we assume without loss of generality that the constant terms in their Taylor expansions are all equal to $1$). Note that $\nu_i(t)$, $i\in\{0,1,2\}$ vanish for $t<0$.   Then, by \cite[Chapter~VI, section~10, formula (1)]{Widder} the convolution  of the measures of this special form is given by
\begin{multline*}
\int\limits_{[0,t)}(\nu_0(u)+\delta_0)du=\int\limits_{[0,t)}\left(\int_{0}^{t-u}\nu_1(x)dx+\mathbf{1}_{[u\le{t}]}\right)\left(\nu_2(u)+\delta_0\right)du
\\
=\int\limits_{[0,t)}\nu_2(u)\left(\int_{0}^{t-u}\nu_1(x)dx\right)du+\int\limits_{[0,t)}\nu_2(u)du+ \int\limits_{[0,t)}\left(\int_{0}^{t-u}\nu_1(x)dx+\mathbf{1}_{[u\le{t}]}\right)\delta_0du
\\
=\int_{0}^{t}\nu_2(u)\left(\int_{0}^{t-u}\nu_1(x)dx\right)du+\int_{0}^{t}\nu_2(u)du+ \int_{0}^{t}\nu_1(x)dx+\mathbf{1}_{[0\le{t}]}.
\end{multline*}
Taking derivatives on both sides and using $\partial_t\mathbf{1}_{[0\le{t}]}=\delta_0$ and
$$
\partial_t\int_{0}^{t}\nu_2(u)\left(\int_{0}^{t-u}\nu_1(x)dx\right)du=\int_{0}^{t}\nu_2(t-u)\nu_1(u)du,
$$
we finally obtain:
\begin{equation}\label{eq:convolution}
\nu_0(t)=\int_{0}^{t}\nu_2(t-u)\nu_1(u)du+\nu_1(t)+\nu_2(t).
\end{equation}
Combined with the Laplace convolution theorem \cite[Theorem~16a]{Widder}  these observations can be summarized as follows.
\begin{lemma}\label{lm:convolution}
Suppose $x\to\phi_i(x)$, $\phi_i(0)=1$, are absolutely monotonic for $i\in\{0,1,2\}$ and are related by $\phi_0(x)=\phi_1(x)\phi_2(x)$.  If $\phi_i(1/x)=\int_{[0,\infty)}e^{-xt}(\nu_i(t)+\delta_0)dt$ for $i\in\{0,1,2\}$ in accordance with Theorem~\ref{th:twoCMfunctions}, then the representing functions $\nu_i$, $i\in\{0,1,2\}$ are related by \emph{(\ref{eq:convolution})}.
\end{lemma}

The following examples illustrate the application of this lemma to derive several new integral evaluations for products of hypergeometric functions.
We believe them to be new. They are also inaccessible for \textit{Mathematica}.  These examples may be viewed as a complement to \cite[section 4.24.5]{BrychkovBook}, containing a number of similar but still different integrals.

\noindent\textbf{Example~1.} Write the identity $(1-x)^{-a}(1-x)^{-b}=(1-x)^{-a-b}$ in the form
$$
{}_1F_{0}\left(\left.\!\!\begin{array}{c}a\\-\end{array}\right|1/x\!\right){}_1F_{0}\left(\left.\!\!\begin{array}{c}b\\-\end{array}\right|1/x\!\right)
={}_1F_{0}\left(\left.\!\!\begin{array}{c}a+b\\-\end{array}\right|1/x\!\right).
$$
Applying (\ref{eq:pFq1-over-x}) to each term above we get explicit expressions for the representing measures $\nu_i$, $i=0,1,2$. Substituting them into
(\ref{eq:convolution}) yields:
\begin{multline*}
ab\!\!\int_{0}^{t}\!\!{}_1F_{1}\!\left(\left.\!\!\!\begin{array}{c}a+1\\2\end{array}\!\!\right|t-u\!\right)
{}_1F_{1}\!\left(\left.\!\!\!\begin{array}{c}b+1\\2\end{array}\!\!\right|u\!\right)du
\\
=(a+b){}_1F_{1}\!\left(\left.\!\!\!\!\begin{array}{c}a+b\\2\end{array}\!\!\right|t\!\right)
-a{}_1F_{1}\!\left(\left.\!\!\!\!\begin{array}{c}a+1\\2\end{array}\!\!\right|t\!\right)
-b{}_1F_{1}\!\left(\left.\!\!\!\!\begin{array}{c}b+1\\2\end{array}\!\!\right|t\!\right).
\end{multline*}

\noindent\textbf{Example~2.} Write Euler transformation \cite[Theorem~2.2.5]{AAR} in the form
$$
{}_1F_{0}\left(\left.\!\!\begin{array}{c}a+b-c\\-\end{array}\right|1/x\!\right){}_2F_{1}\left(\left.\!\!\begin{array}{c}c-a,c-b\\c\end{array}\right|1/x\!\right)
={}_2F_{1}\left(\left.\!\!\begin{array}{c}a,b\\c\end{array}\right|1/x\!\right).
$$
Applying (\ref{eq:pFq1-over-x}) to each term above we get explicit expressions for the representing measures $\nu_i$, $i=0,1,2$. Substituting them into
(\ref{eq:convolution}) yields:
\begin{multline*}
(a+b-c)\!\!\int_{0}^{t}\!\!{}_1F_{1}\!\left(\left.\!\!\!\!\begin{array}{c}a+b-c+1\\2\end{array}\!\!\right|t-u\!\right)
{}_2F_{2}\!\left(\left.\!\!\!\!\begin{array}{c}c-a+1,c-b+1\\c+1,2\end{array}\!\!\right|u\!\right)du
\\
=\frac{ab}{(c-a)(c-b)}{}_2F_{2}\!\left(\left.\!\!\!\!\begin{array}{c}a+1,b+1\\c+1,2\end{array}\!\!\right|t\!\right)
-\frac{c(a+b-c)}{(c-a)(c-b)}{}_1F_{1}\!\left(\left.\!\!\!\!\begin{array}{c}a+b-c+1\\2\end{array}\!\!\right|t\!\right)
\\
-{}_2F_{2}\!\left(\left.\!\!\!\!\begin{array}{c}c-a+1,c-b+1\\c+1,2\end{array}\!\!\right|t\!\right).
\end{multline*}

\medskip

\noindent\textbf{Example~3.}  Take the celebrated Clausen's identity \cite[(11)]{Grinshpan} in the form
$$
\left[{}_2F_{1}\left(\left.\!\!\begin{array}{c}a-1,b-1\\a+b-3/2\end{array}\right|1/x\!\right)\right]^2
={}_3F_{2}\left(\left.\!\!\begin{array}{c}2a-2,2b-2,a+b-2\\a+b-3/2,2a+2b-4\end{array}\right|1/x\!\right).
$$
Writing the representing measure for each function in this formula by (\ref{eq:pFq1-over-x}) and applying (\ref{eq:convolution}) we obtain
after elementary manipulations:
\begin{multline*}
\frac{(a-1)(b-1)}{2a+2b-3}\int_{0}^{t}{}_2F_{2}\left(\left.\!\!\begin{array}{c}a,b\\a+b-1/2,2\end{array}\right|t-u\!\right)
{}_2F_{2}\left(\left.\!\!\begin{array}{c}a,b\\a+b-1/2,2\end{array}\right|u\!\right)du
\\
={}_3F_{3}\left(\left.\!\!\begin{array}{c}2a-1,2b-1,a+b-1\\a+b-1/2,2a+2b-3,2\end{array}\right|t\!\right)
-{}_2F_{2}\left(\left.\!\!\begin{array}{c}a,b\\a+b-1/2,2\end{array}\right|t\!\right).
\end{multline*}

\noindent\textbf{Example~4.}  Next, we take the product formula for the Bessel functions in the form \cite[(18)]{Grinshpan}
$$
{}_0F_{1}\left(\left.\!\!\begin{array}{c}-\\a\end{array}\right|1/x\!\right){}_0F_{1}\left(\left.\!\!\begin{array}{c}-\\b\end{array}\right|1/x\!\right)
={}_2F_{3}\left(\left.\!\!\begin{array}{c}(a+b)/2,(a+b-1)/2\\a,b,a+b-1\end{array}\right|4/x\!\right).
$$
As the argument on the right hand side is $4/x$, we need a slightly modified form of (\ref{eq:pFq1-over-x}) given by
$$
{_{p}F_q}\left(\a;\b;\alpha/x\right)=\int\limits_{[0,\infty)}e^{-xt}\left[\alpha\frac{(\a)}{(\b)}
{_{p}F_{q+1}}\left(\left.\!\!\begin{array}{c}\a+1\\\b+1,2\end{array}\right|\alpha{t}\!\right)+\delta_0\right]dt
$$
and obtained by the obvious change of variable.  This formula combined with (\ref{eq:convolution}) yields:
\begin{multline*}
\frac{1}{ab}\int_{0}^{t}{}_0F_{2}\left(\left.\!\!\begin{array}{c}-\\a+1,2\end{array}\right|t-u\!\right)
{}_0F_{2}\left(\left.\!\!\begin{array}{c}-\\b+1,2\end{array}\right|u\!\right)du
\\
=\frac{(a+b)}{ab}{}_2F_{4}\!\left(\left.\!\!\begin{array}{c}(a+b)/2+1,(a+b-1)/2+1\\a+1,b+1,a+b,2\end{array}\!\right|4t\!\right)
-\frac{1}{a}{}_0F_{2}\!\left(\left.\!\!\begin{array}{c}-\\a+1,2\end{array}\!\right|t\!\right)
-\frac{1}{b}{}_0F_{2}\!\left(\left.\!\!\begin{array}{c}-\\b+1,2\end{array}\!\right|t\!\right).
\end{multline*}

\noindent\textbf{Example~5.}  Using the same algorithm, Orr's identities \cite[(14)-(16)]{Grinshpan} lead to the next three integral evaluations:
\begin{multline*}
ab\int_{0}^{t}{}_2F_{2}\left(\left.\!\!\begin{array}{c}a+1,b+1\\a+b+1/2,2\end{array}\right|t-u\!\right)
{}_2F_{2}\left(\left.\!\!\begin{array}{c}a+1,b+1\\a+b+3/2,2\end{array}\right|u\!\right)du
\\
=2(a+b){}_3F_{3}\!\left(\left.\!\!\begin{array}{c}2a+1,2b+1,a+b+1\\a+b+3/2,2a+2b,2\end{array}\!\right|t\!\right)
-(a+b+1/2){}_2F_{2}\!\left(\left.\!\!\begin{array}{c}a+1,b+1\\a+b+1/2,2\end{array}\!\right|t\!\right)
\\
-(a+b-1/2){}_2F_{2}\!\left(\left.\!\!\begin{array}{c}a+1,b+1\\a+b+3/2,2\end{array}\!\right|t\!\right),
\end{multline*}
\begin{multline*}
\frac{ab(b-1)}{(a+b-1/2)}\int_{0}^{t}{}_2F_{2}\left(\left.\!\!\begin{array}{c}a+1,b+1\\a+b+1/2,2\end{array}\right|t-u\!\right)
{}_2F_{2}\left(\left.\!\!\begin{array}{c}a+1,b\\a+b+1/2,2\end{array}\right|u\!\right)du
\\
=(2b-1){}_3F_{3}\!\left(\left.\!\!\begin{array}{c}2a+1,2b,a+b\\a+b+1/2,2a+2b-1,2\end{array}\!\right|t\!\right)
-b{}_2F_{2}\!\left(\left.\!\!\begin{array}{c}a+1,b+1\\a+b+1/2,2\end{array}\!\right|t\!\right)
\\
-(b-1){}_2F_{2}\!\left(\left.\!\!\begin{array}{c}a+1,b\\a+b+1/2,2\end{array}\!\right|t\!\right)
\end{multline*}
and
\begin{multline*}
\frac{ab(1/2-a)(1/2-b)}{(a+b+1/2)(3/2-a-b)}\int_{0}^{t}{}_2F_{2}\left(\left.\!\!\begin{array}{c}a+1,b+1\\a+b+3/2,2\end{array}\right|t-u\!\right)
{}_2F_{2}\left(\left.\!\!\begin{array}{c}3/2-a,3/2-b\\5/2-a-b,2\end{array}\right|u\!\right)du
\\
=\frac{(a-b+1/2)(b-a+1/2)}{2(a+b+1/2)(3/2-a-b)}{}_3F_{3}\!\left(\left.\!\!\begin{array}{c}3/2,a-b+3/2,b-a+3/2\\a+b+3/2,5/2-a-b,2\end{array}\!\right|t\!\right)
\\
-\frac{ab}{(a+b+1/2)}{}_2F_{2}\!\left(\left.\!\!\begin{array}{c}a+1,b+1\\a+b+3/2,2\end{array}\!\right|t\!\right)
-\frac{(1/2-a)(1/2-b)}{(3/2-a-b)}{}_2F_{2}\!\left(\left.\!\!\begin{array}{c}3/2-1,3/2-b\\5/2-a-b,2\end{array}\!\right|t\!\right).
\end{multline*}

\section{Stieltjes transform representation and its consequences}

Recall that a function $f:(0,\infty)\to\R$ is called a Stieltjes function, if it is of the form
$$
f(x)=c+\int\nolimits_{[0,\infty)}\frac{\mu(dt)}{x+t},
$$
where $c$ is a nonnegative constant and $\mu$ is a positive measure on $[0, \infty)$ such that the above integral converges for any $x>0$, see \cite[Chapter~VIII]{Widder} and \cite[p.294]{Koumandos}. In \cite[Theorem~8]{KPJMAA} we found sufficient conditions on parameters in order that the function ${_{p+1}F_p}(\a;\b;-x)$ be of this class.  This allows the use of the Stieltjes inversion formula \cite[Theorem 7a]{Widder} to calculate the jump of the generalized hypergeometric function ${_{p+1}F_p}(z)$ over the branch cut $[1,\infty)$.  We can then employ an analytic continuation formula for the generalized hypergeometric function to get rid of the restrictions on parameters.   Notwithstanding the simplicity of the calculation that follows, we are unaware of any references for formula (\ref{eq:jump}).
\begin{theorem}\label{th:jump}
Suppose $x>1$ and  $\a,\b$ are real vectors.  Then the following identities hold true
\begin{equation}\label{eq:jump}
{_{p+1}F_p}\!\left(\!\!\begin{array}{l}\a\\\b\end{array}\vline\,\,x+i0\!\right)
-{_{p+1}F_p}\!\left(\!\!\begin{array}{l}\a\\\b\end{array}\vline\,\,x-i0\!\right)
=2\pi{i}\frac{\Gamma(\b)}{\Gamma(\a)}G^{p+1,0}_{p+1,p+1}\left(\frac{1}{x}\,\vline\begin{array}{c}1,\b\\\a\end{array}\!\!\right),
\end{equation}

\begin{equation}\label{eq:avver}
\frac{{_{p+1}F_p}\left(\a;\b;x+i0\right)+{_{p+1}F_p}\left(\a;\b;x-i0\right)}{2}
\\
=-\frac{\pi\Gamma(\b)}{\sqrt{x}\Gamma(\a)}
G^{p+1,1}_{p+2,p+2}\left(\frac{1}{x}\,\vline\begin{array}{c}1/2,1,\b-1/2\\\a-1/2,1\end{array}\!\!\right).
\end{equation}
\end{theorem}

\textbf{Proof.} According to \cite[Theorem~8]{KPJMAA} the next formula holds
$$
{_{p+1}F_p}\!\!\left(\!\!\begin{array}{l}\a\\\b\end{array}\vline\,\,-z\!\right)=\int_{1}^{\infty}\frac{\rho(x)dx}{x+z},
$$
for $0<a_1\le1$, $\b\prec^W\a_{[1]}$,  where
$$
\rho(x)=\frac{\Gamma(\b)}{\Gamma(\a)}G^{p+1,0}_{p+1,p+1}\left(\frac{1}{x}\,\vline\begin{array}{c}1,\b\\\a\end{array}\!\!\right).
$$
According to the Stieltjes inversion formula \cite[Theorem 7a]{Widder} or \cite[(12)]{Schwarz}
$$
\rho(x)=\frac{1}{2\pi{i}}\left\{{_{p+1}F_p}\!\!\left(\!\!\begin{array}{l}\a\\\b\end{array}\vline\,\,x+i0\!\right)
-{_{p+1}F_p}\!\!\left(\!\!\begin{array}{l}\a\\\b\end{array}\vline\,\,x-i0\!\right)\right\}
$$
and the claim follows under the above restrictions on parameters. To get rid of these restrictions we can apply \cite[formula (16.8.6)]{NIST} reading
\begin{equation}\label{eq:hypergeom}
{_{p+1}F_p}\!\!\left(\!\!\begin{array}{l}\a\\\b\end{array}\vline\,\,x\!\right)
=\sum\limits_{j=1}^{p+1}\frac{\Gamma(\a_{[j]}-a_j)\Gamma(\b)}{\Gamma(\a_{[j]})\Gamma(\b-a_j)}
(-x)^{-a_j}{_{p+1}F_p}\!\!\left(\!\!\begin{array}{l}a_j,1-\b+a_j\\1-\a_{[j]}+a_j\end{array}\vline\,\,\frac{1}{x}\!\right).
\end{equation}
Since  $(-x-i0)^{-a_j}-(-x+i0)^{-a_j}={2i}x^{-a_j}\sin(\pi{a_j})$ for $x>1$,  we get
\begin{multline*}
{_{p+1}F_p}\!\!\left(\!\!\begin{array}{l}\a\\\b\end{array}\vline\,\,x+i0\!\right)
-{_{p+1}F_p}\!\!\left(\!\!\begin{array}{l}\a\\\b\end{array}\vline\,\,x-i0\!\right)
\\
=2i\sum\limits_{j=1}^{p+1}\frac{\Gamma(\a_{[j]}-a_j)\Gamma(\b)}{\Gamma(\a_{[j]})\Gamma(\b-a_j)}x^{-a_j}\sin(\pi{a_j})
{_{p+1}F_p}\!\left(\!\!\begin{array}{l}a_j,1-\b+a_j\\1-\a_{[j]}+a_j\end{array}\!\!\vline\,\,\frac{1}{x}\right).
\end{multline*}
The identity
$$
2i\sum\limits_{j=1}^{p+1}\frac{\Gamma(\a_{[j]}-a_j)\Gamma(\b)}{\Gamma(\a_{[j]})\Gamma(\b-a_j)}x^{-a_j}\sin(\pi{a_j})
{_{p+1}F_p}\!\left(\!\!\begin{array}{l}a_j,1-\b+a_j\\1-\a_{[j]}+a_j\end{array}\!\!\vline\,\,\frac{1}{x}\right)
=2\pi{i}\frac{\Gamma(\b)}{\Gamma(\a)}G^{p+1,0}_{p+1,p+1}\!\left(\frac{1}{x}\,\vline\begin{array}{c}1,\b\\\a\end{array}\!\!\right)
$$
following from  \cite[(2.4)]{KPSigma} by Euler's reflection formula completes the proof
of formula (\ref{eq:jump}). To demonstrate (\ref{eq:avver}) we apply
$$
(-x-i0)^{-a_j}+(-x+i0)^{-a_j}={2}x^{-a_j}\cos\pi{a_j}=-\frac{\pi}{\Gamma(a_j-1/2)\Gamma(3/2-a_j)}
$$
together with formula (\ref{eq:hypergeom}) to get (assuming  $x>1$):
\begin{multline*}
\frac{{_{p+1}F_p}\left(\a;\b;x+i0\right)+{_{p+1}F_p}\left(\a;\b;x-i0\right)}{2}
\\
=-\frac{\pi}{\sqrt{x}}\sum\limits_{j=1}^{p+1}\frac{\Gamma(\a_{[j]}-a_j)\Gamma(\b)(1/x)^{a_j-1/2}}{\Gamma(\a_{[j]})\Gamma(\b-a_j)\Gamma(a_j-1/2)\Gamma(3/2-a_j)}
{_{p+2}F_{p+1}}\!\left(\!\begin{array}{l}1-\b+a_j,a_j,a_j-1/2\\1-\a_{[j]}+a_j,a_j-1/2\end{array}\!\!\vline\,\,\frac{1}{x}\right).
\end{multline*}
To verify that this expression is equal to the right hand side of (\ref{eq:avver}) it remains to apply \cite[16.17.2]{NIST} or \cite[8.2.2.3]{PBM3}.
$\hfill\square$

As a corollary we recover \cite[(15.2.3)]{NIST}:
\begin{corollary}\label{cr:2F1jump}
Suppose $x>1.$  Then the following identity holds
\begin{multline}\label{eq:2F1jump}
{}_2F_1\!\left(\!\!\begin{array}{l} a,b\\c\end{array}\vline\,\,x+i0\!\right)-{}_2F_1\!\left(\!\!\begin{array}{l}a,b\\c\end{array}\vline\,\,x-i0\!\right)
\\
=\frac{2\pi{i}\Gamma(c)}{\Gamma(a)\Gamma(b)\Gamma(1+c-a-b)}(x-1)^{c-a-b}{}_2F_1\!\left(\!\!\begin{array}{l}c-a,c-b\\c-a-b+1\end{array}\vline\,\,1-x\!\right).
\end{multline}
\end{corollary}
\textbf{Proof}.  Indeed, writing (\ref{eq:jump}) for $p=1$ gives
$$
{_{2}F_1}\!\left(\!\!\begin{array}{l} a,b\\c\end{array}\vline\,\,x+i0\!\right)-{_{2}F_1}\!\left(\!\!\begin{array}{l}a, b\\ c\end{array}\vline\,\,x-i0\!\right)
=2\pi{i}\frac{\Gamma(c)}{\Gamma(a)\Gamma(b)}G^{2,0}_{2,2}\left(\frac{1}{x}\,\vline\begin{array}{c}1,c\\a,b\end{array}\!\!\right).
$$
Further expressing $G^{2,0}_{2,2}\left(1/x\,\vline\begin{array}{c}1,c\\a,b\end{array}\!\!\right)$ via ${}_2F_1$ according to \cite[p.4]{KPSigma} and using Pfaff's
transformation \cite[(15.8.1)]{NIST}  we arrive at (\ref{eq:2F1jump}). $\hfill\square$

\textbf{Remark.} Formulas (\ref{eq:jump}) and (\ref{eq:avver}) are consistent with the boundary values of the generalized hypergeometric function
given in \cite{Wolfram1}. In fact, these reference can be used to give an alternative derivation of (\ref{eq:jump}) and (\ref{eq:avver}). Indeed, we haveÑ
\begin{multline*}
\frac{\Gamma(\a)}{\Gamma(\b)}\left({_{p+1}F_p}\!\left(\!\!\begin{array}{l}\a\\\b\end{array}\vline\,\,x+i0\!\right)
-{_{p+1}F_p}\!\left(\!\!\begin{array}{l}\a\\\b\end{array}\vline\,\,x-i0\!\right)\right)
\\
=G^{p+1,1}_{p+1,p+1}\left(e^{i\pi}\frac{1}{x}\,\vline\begin{array}{c}1,\b\\\a\end{array}\!\!\right)
-G^{p+1,1}_{p+1,p+1}\left(e^{-i\pi}\frac{1}{x}\,\vline\begin{array}{c}1,\b\\\a\end{array}\!\!\right)
\\
=\frac{-1}{2\pi{i}}\int\limits_{\L}\frac{e^{i\pi{s}}\Gamma(-s)\Gamma(s+\a)}{\Gamma(s+\b)}\left(\frac{1}{x}\right)^{-s}ds
+\frac{1}{2\pi{i}}\int\limits_{\L}\frac{e^{-i\pi{s}}\Gamma(-s)\Gamma(s+\a)}{\Gamma(s+\b)}\left(\frac{1}{x}\right)^{-s}ds
\\
=-\frac{1}{\pi}\int\limits_{\L}\frac{e^{i\pi{s}}-e^{-i\pi{s}}}{2i}\frac{\Gamma(-s)\Gamma(s+\a)}{\Gamma(s+\b)}\left(\frac{1}{x}\right)^{-s}ds
=-\frac{1}{\pi}\int\limits_{\L}\frac{\sin(\pi{s})\Gamma(-s)\Gamma(s+\a)}{\Gamma(s+\b)}\left(\frac{1}{x}\right)^{-s}ds
\\
=\int\limits_{\L}\frac{\Gamma(s+\a)}{\Gamma(s+1)\Gamma(s+\b)}\left(\frac{1}{x}\right)^{-s}ds
=2\pi{i}G^{p+1,0}_{p+1,p+1}\left(\frac{1}{x}\,\vline\begin{array}{c}1,\b\\\a\end{array}\!\!\right).
\end{multline*}
In a similar fashion, for the mean value on the cut, we obtain:
\begin{multline*}
\frac{{_{p+1}F_p}\!\left(\!\!\begin{array}{l}\a\\\b\end{array}\vline\,\,x+i0\!\right)+{_{p+1}F_p}\!\left(\!\!\begin{array}{l}\a\\\b\end{array}\vline\,\,x-i0\!\right)}{2}
\\
=\frac{\Gamma(\b)}{2\Gamma(\a)}\left(G^{p+1,1}_{p+1,p+1}\left(e^{i\pi}\frac{1}{x}\,\vline\begin{array}{c}1,\b\\\a\end{array}\!\!\right)+
G^{p+1,1}_{p+1,p+1}\left(e^{-i\pi}\frac{1}{x}\,\vline\begin{array}{c}1,\b\\\a\end{array}\!\!\right)\right)
\\
=\!\frac{\Gamma(\b)}{2\pi{i}\Gamma(\a)}\!\!\int\limits_{\L}\frac{e^{i\pi{s}}+e^{-i\pi{s}}}{2}\frac{\Gamma(-s)\Gamma(s+\a)}
{\Gamma(s+\b)}\left(\frac{1}{x}\right)^{-s}\!\!\!ds
\!=\!\frac{\Gamma(\b)}{2\pi{i}\Gamma(\a)}\!\!\int\limits_{\L}\!\!\cos{\pi{s}}\frac{\Gamma(-s)\Gamma(s+\a)}{\Gamma(s+\b)}\left(\frac{1}{x}\right)^{-s}\!\!\!ds
\\
=[s=t-1/2]=\frac{\Gamma(\b)\pi}{2\pi{i}\sqrt{x}\Gamma(\a)}\int\limits_{\L'}\frac{\Gamma(-t)\sin{\pi{t}}}{\pi}\frac{\Gamma(-t+1/2)\Gamma(t+\a-1/2)}
{\Gamma(-t)\Gamma(t+\b-1/2)}\left(\frac{1}{x}\right)^{-t}dt
\\
=\frac{\Gamma(\b)\pi}{2\pi{i}\sqrt{x}\Gamma(\a)}\int\limits_{\L'}\frac{\Gamma(-t)\sin{\pi{t}}}{\pi}\frac{\Gamma(-t+1/2)\Gamma(t+\a-1/2)}
{\Gamma(-t)\Gamma(t+\b-1/2)}\left(\frac{1}{x}\right)^{-t}dt
\\
=-\frac{\Gamma(\b)\pi}{2\pi{i}\sqrt{x}\Gamma(\a)}\int\limits_{\L'}\frac{\Gamma(-t+1/2)\Gamma(t+\a-1/2)}{\Gamma(1+t)\Gamma(-t)\Gamma(t+\b-1/2)}
\left(\frac{1}{x}\right)^{-t}dt
\\
=-\frac{\Gamma(\b)\pi}{\sqrt{x}\Gamma(\a)}G^{p+1,1}_{p+2,p+2}\left(\frac{1}{x}\,\vline\begin{array}{c}1/2,1,\b-1/2\\\a-1/2,1\end{array}\!\!\right).
\end{multline*}
Here $\L'$ denotes the contour indented by $1/2$ with respect to the original contour $\L$.  Wolfram function site \cite{Wolfram1}  defines the value of the generalized hypergeometric function on the cut to be equal to ${_{p+1}F_p}(\a;\b; x-i0)$ by continuity from below.
In our opinion, it makes more sense to define it as the average value given by the right hand side of (\ref{eq:avver}).

\medskip

In \cite[Corollary~1]{KPJMAA} we observed that the function $x\mapsto{_{p+1}F_p}(\sigma,\a;\b;-z^{1/\sigma})$
is a Stiletjes function for each $\sigma>0$ as long as $v_{\a,\b}(t)\ge0$ on $[0,1]$ and derived the representation
\begin{equation}\label{eq:q1Fqnewrepr}
{_{p+1}F_p}(\sigma,\a;\b;-z)=\int\limits_{0}^{\infty}\frac{\varphi(y)dy}{y^{\sigma}+z^{\sigma}},
\end{equation}
valid for $\sigma\ge2$ and $|\arg(z)|<\pi/\sigma$.  The density $\varphi$ is given by
\begin{equation}\label{eq:varphigeneral}
\varphi(y)=\frac{\sigma{y^{\sigma-1}}\Gamma(\b)}{\pi\Gamma(\a)}
\int\limits_{0}^{1}\frac{\sin\left\{\sigma\arctan\left(\frac{ty\sin(\pi/\sigma)}{1+ty\cos(\pi/\sigma)}\right)\right\}}
{t\left(1+2ty\cos(\pi/\sigma)+t^2y^2\right)^{\sigma/2}}G^{p,0}_{p,p}\left(t\,\,\vline\begin{array}{c}\!\b\\\!\a\end{array}\!\!\right)dt.
\end{equation}
For $\sigma=2$ formula (\ref{eq:varphigeneral}) simplifies to
\begin{equation}\label{eq:varphi2}
\varphi(y)=\frac{4\Gamma(\b)}{\pi\Gamma(\a)}\int\limits_{0}^{1}\frac{y^2}{\left(1+t^2y^2\right)^{2}}
G^{p,0}_{p,p}\left(t\,\vline\begin{array}{c}\!\b\\\!\a\end{array}\!\!\right)\!dt
\end{equation}
by virtue of the identity $\sin(2\arctan(s))=2s/(1+s^2)$. Below we simplify the expression for $\varphi$ further.
\begin{theorem}\label{th:sigmaequal2}
Suppose $\Re(\a)>0$, then
\begin{equation}\label{eq:sigmaequal2}
{_{p+1}F_{p}}\left(\left.\!\!\begin{array}{c}\a\\\b\end{array}\right|-z\!\right)
=\frac{2(\a)}{\pi(\b)}\int\nolimits_{0}^{\infty}\frac{t^2}{z^2+t^2}{_{2p+2}F_{2p+1}}\!\!\left(\left.\!\!\begin{array}{c}(\a+1)/2,(\a+2)/2\\(\b+1)/2,(\b+2)/2,3/2\end{array}\right|-t^2\!\right)\!dt
\end{equation}
for $|\arg(z)|<\pi/2$.
\end{theorem}
\textbf{Proof.}  Set $\sigma=2$ in (\ref{eq:q1Fqnewrepr}) and calculate using the binomial expansion in (\ref{eq:varphi2}):
\begin{multline*}
\varphi(y)=\frac{4\Gamma(\b)}{\pi\Gamma(\a)}\int\nolimits_{0}^{1}\frac{y^2}{\left(1+t^2y^2\right)^{2}}
G^{p,0}_{p,p}\left(t\,\vline\begin{array}{c}\!\b\\\!\a\end{array}\!\!\right)\!dt
=\frac{4y^2\Gamma(\b)}{\pi\Gamma(\a)}\!\sum\limits_{k=0}^{\infty}\frac{(2)_k(-y^2)^{k}}{k!}
\!\int\nolimits_{0}^{1}\!G^{p,0}_{p,p}\left(t\,\vline\begin{array}{c}\!\b\\\!\a\end{array}\!\!\right)t^{2k}\!dt
\\
=\frac{4y^2\Gamma(\b)}{\pi\Gamma(\a)}\!\sum\limits_{k=0}^{\infty}\frac{(2)_k(-y^2)^{k}}{k!}
\frac{\Gamma(\a+2k+1)}{\Gamma(\b+2k+1)}
\\
=\frac{4y^2\Gamma(\b)}{\pi\Gamma(\a)}\!\sum\limits_{k=0}^{\infty}\frac{(2)_k}{k!}
\frac{2^{\a+2k}\Gamma(\a/2+k+1/2)\Gamma(\a/2+k+1)}{2^{\b+2k}\Gamma(\b/2+k+1/2)\Gamma(\b/2+k+1)}(-y^2)^{k}
\\
=\frac{4y^22^{\a}\Gamma(\b)\Gamma(\a/2+1/2)\Gamma(\a/2+1)}{\pi2^{\b}\Gamma(\a)\Gamma(\b/2+1/2)\Gamma(\b/2+1)}\!\sum\limits_{k=0}^{\infty}\frac{(2)_k}{k!}
\frac{(\a/2+1/2)_{k}(\a/2+1)_{k}}{(\b/2+1/2)_{k}(\b/2+1)_{k}}(-y^2)^{k}
\\
=\frac{4y^2(\a)}{\pi(\b)}{_{2p+1}F_{2p}}\left(\left.\!\!\begin{array}{c}2,(\a+1)/2,(\a+2)/2\\(\b+1)/2,(\b+2)/2\end{array}\right|-y^2\!\right).
\end{multline*}
Substituting into (\ref{eq:q1Fqnewrepr}) we get:
\begin{multline*}
{_{p+1}F_p}(\a;\b;-z)={_{p+2}F_{p+1}}(2,\a;2,\b;-z)
\\
=\frac{4(\a)}{\pi2(\b)}\int\limits_{0}^{\infty}{_{2p+3}F_{2p+2}}\left(\left.\!\!\begin{array}{c}2,(\a+1)/2,(\a+2)/2\\(\b+1)/2,(\b+2)/2,(2+1)/2,(2+2)/2\end{array}\right|-y^2\!\right)\frac{y^2dy}{y^{2}+z^{2}}
\\
=\frac{2(\a)}{\pi(\b)}\int\limits_{0}^{\infty}{_{2p+2}F_{2p+1}}\left(\left.\!\!\begin{array}{c}(\a+1)/2,(\a+2)/2\\(\b+1)/2,(\b+2)/2,3/2\end{array}\right|-y^2\!\right)\frac{y^2dy}{y^{2}+z^{2}}.
\end{multline*}
The convergence condition $\Re(\a)>0$ is verified by (\ref{eq:hypergeom}). $\hfill\square$

\textbf{Remark.} Formula (\ref{eq:sigmaequal2}) is not contained in \cite{PBM3} and \textit{Mathematica} is unable to evaluate the integral on the right hand side.

We can deduce another simplification of (\ref{eq:varphigeneral}) for $\sigma=3$.
\begin{theorem}\label{th:sigmaequal3}
Suppose $\Re(\a)>0$, then for $|\arg(z)|<\pi/3$
\begin{equation}\label{eq:sigmaequal3}
{_{p+1}F_{p}}\left(\left.\!\!\begin{array}{c}\a\\\b\end{array}\right|-z\!\right)
=\frac{9\sqrt{3}}{2\pi}\int\nolimits_{0}^{\infty}\frac{\psi_1(y)-2\psi_2(y)+2\psi_4(y)-\psi_5(y)}{z^3+y^3}dy,
\end{equation}
where
$$
\psi_m(y)=\frac{y^{m+2}(\a)_m}{(3)_m(\b)_m}
{_{3p+4}F_{3p+3}}\left(\left.\!\!\begin{array}{c}3,\Delta(\a+m,3)\\\Delta(\b+m,3),\Delta(m+3,3)\end{array}\right|y^3\!\right),
$$
and
$\Delta(\mathbf{c},k)=(\mathbf{c}/k,(\mathbf{c}+1)/k,\ldots,(\mathbf{c}+k-1)/k)$.
The values of the hypergeometric function ${}_{p+1}F_{p}(t)$ for $t\in[1,\infty)$ can be taken on either side of the branch cut \emph{(}but consistently for all $\psi_m$\emph{)} or as the average value according to formula \emph{(\ref{eq:avver})}. In either case, the numerator of the integrand in \emph{(\ref{eq:sigmaequal3})} takes real values.
\end{theorem}
\textbf{Proof.}  First, we treat the case of ${_{p+1}F_{p}}(3,\a;\b;-z)$, where $\a$ and $\b$ both have $p$ components.  Taking $\sigma=3$ in (\ref{eq:q1Fqnewrepr}) we want to calculate $\varphi(y)$ in (\ref{eq:varphigeneral}). Using the formulas \cite[(18.5.2)]{NIST}
$$
\sin(n\theta)=\sin(\theta)U_{n-1}(\cos(\theta))~~\text{and}~~\cos(\arctan(x))=(1+x^2)^{-1/2},
$$
where $U_{n}$ denotes the Chebyshev polynomial of the second kind, and \cite[(18.5.15)]{NIST} we get
$$
\sin(3\arctan(x))=\frac{x(3-x^2)}{(1+x^2)^{3/2}}.
$$
Tedious but elementary calculations then yield
\begin{equation}\label{eq:sin3arctan}
\frac{\sin\left(3\arctan\left(\frac{ty\sin(\pi/3)}{1+ty\cos(\pi/3)}\right)\right)}{(1+2ty\cos(\pi/3)+t^2y^2)^{3/2}}
=\frac{3\sqrt{3}ty(1+ty)}{2(1+ty+t^2y^2)^3}=\frac{3\sqrt{3}ty(1-2ty+2t^3y^3-t^4y^4)}{2(1-t^3y^3)^3}.
\end{equation}
Applying the binomial expansion to $(1-t^3y^3)^{-3}$ and integrating termwise we obtain:
\begin{multline*}
\varphi(y)=\frac{9\sqrt{3}y^2\Gamma(\b)}{2\pi\Gamma(\a)}\int\limits_{0}^{1}ty(1-2ty+2t^3y^3-t^4y^4)\sum\limits_{k=0}^{\infty}
\frac{(3)_k}{k!}(ty)^{3k}G^{p,0}_{p,p}\left(t\,\,\vline\begin{array}{c}\!\b-1\\\!\a-1\end{array}\!\!\right)dt
\\
=\frac{9\sqrt{3}y^2\Gamma(\b)}{2\pi\Gamma(\a)}(\varphi_1(y)-2\varphi_2(y)+2\varphi_4(y)-\varphi_5(y)),
\end{multline*}
where
\begin{multline*}
\varphi_m(y)=\int\limits_{0}^{1}\sum\limits_{k=0}^{\infty}\frac{(3)_k}{k!}(ty)^{3k+m}G^{p,0}_{p,p}\left(t\,\,\vline\begin{array}{c}\!\b-1\\\!\a-1\end{array}\!\!\right)dt
=\sum\limits_{k=0}^{\infty}\frac{(3)_k}{k!}y^{3k+m}\int\limits_{0}^{1}t^{3k+m-1}G^{p,0}_{p,p}\left(t\,\,\vline\begin{array}{c}\!\b\\\!\a\end{array}\!\!\right)dt
\\
=\sum\limits_{k=0}^{\infty}\frac{(3)_k}{k!}y^{3k+m}\int\limits_{0}^{1}t^{3k+m-1}G^{p,0}_{p,p}\left(t\,\,\vline\begin{array}{c}\!\b\\\!\a\end{array}\!\!\right)dt
=\sum\limits_{k=0}^{\infty}\frac{(3)_k}{k!}y^{3k+m}\frac{\Gamma(\a+3k+m)}{\Gamma(\b+3k+m)}
\\
=3^{\sum_{i=1}^{p}(a_i-b_i)}\sum\limits_{k=0}^{\infty}\frac{(3)_k}{k!}y^{3k+m}
\frac{\Gamma((\a+m)/3+k)\Gamma((\a+m+1)/3+k)\Gamma((\a+m+2)/3+k)}{\Gamma((\b+m)/3+k)\Gamma((\b+m+1)/3+k)\Gamma((\b+m+2)/3+k)}
\\
=y^m\frac{\Gamma(\a+m)}{\Gamma(\b+m)}
{_{3p+1}F_{3p}}\left(\left.\!\!\begin{array}{c}3,(\a+m)/3,(\a+m+1)/3,(\a+m+2)/3\\(\b+m)/3,(\b+m+1)/3,(\b+m+2)/3\end{array}\right|y^3\!\right).
\end{multline*}
According to (\ref{eq:q1Fqnewrepr}) the function $\varphi(y)$ is the density in the representation of ${_{p+1}F_{p}}(3,\a;\b;-z)$. For the general case just write
${_{p+1}F_{p}}(\a;\b;-z)={_{p+2}F_{p+1}}(3,\a;3,\b;-z)$ and apply (\ref{eq:q1Fqnewrepr}) with $\b$ substituted by $(3,\b)$. This leads immediately to formula (\ref{eq:sigmaequal3}).  To verify the reality of the numerator in (\ref{eq:sigmaequal3}) note that we can rewrite $\varphi(y)$ using the first formula in (\ref{eq:sin3arctan}):
$$
\varphi(y)=\frac{9\sqrt{3}y^3\Gamma(\b)}{2\pi\Gamma(\a)}\int\limits_{0}^{1}\frac{1+ty}{(1+ty+t^2y^2)^3}G^{p,0}_{p,p}\left(t\,\,\vline\begin{array}{c}\!\b\\\!\a\end{array}\!\!\right)dt.
$$
The integral here is well-defined and real for all $y>0$.  The convergence condition $\Re(\a)>0$ can be verified by application of (\ref{eq:hypergeom}) to each $\psi_m$ in the integrand.  $\hfill\square$

\textbf{Remark.} The above formula for $\varphi(y)$ permits derivation of another expression for this function as a double hypergeometric series.  Indeed, writing
\begin{multline*}
\frac{1+ty}{(1+ty+t^2y^2)^3}=\frac{1+ty}{(1+ty(1+ty))^3}=\sum\limits_{k=0}^{\infty}(-1)^k\frac{(3)_k}{k!}(ty)^k(1+ty)^{k+1}
\\
=\sum\limits_{k=0}^{\infty}(-1)^k\frac{(3)_k}{k!}\sum\limits_{n=0}^{k+1}\binom{k+1}{n}(ty)^{n+k},
\end{multline*}
we can substitute this into the above integral and integrate termwise to get
\begin{multline*}
\varphi(y)=\frac{9\sqrt{3}y^3\Gamma(\b)}{2\pi\Gamma(\a)}\sum\limits_{k=0}^{\infty}(-1)^k\frac{(3)_k}{k!}y^{n+k}\sum\limits_{n=0}^{k+1}\binom{k+1}{n}
\int\limits_{0}^{1}t^{n+k}G^{p,0}_{p,p}\left(t\,\,\vline\begin{array}{c}\!\b\\\!\a\end{array}\!\!\right)dt
\\
=\frac{9\sqrt{3}y^3\Gamma(\b)}{2\pi\Gamma(\a)}\sum\limits_{k=0}^{\infty}(-1)^k\frac{(3)_k}{k!}y^{n+k}\sum\limits_{n=0}^{k+1}\binom{k+1}{n}
\frac{\Gamma(\a+n+k+1)}{\Gamma(\b+n+k+1)}
\\
=\frac{9\sqrt{3}y^3}{2\pi}
\sum\limits_{j=0}^{\infty}\frac{(\a)_{j+1}}{(\b)_{j+1}}y^{j}\sum\limits_{n+k=j}(-1)^k\frac{(3)_k}{k!}\binom{k+1}{n}
\\
=\frac{9\sqrt{3}y^3}{2\pi}
\sum\limits_{j=0}^{\infty}\frac{(\a)_{j+1}}{(\b)_{j+1}}y^{j}\sum\limits_{k=[(j-1)/2]}^{j}(-1)^k\frac{(3)_k}{k!}\binom{k+1}{j-k}.
\end{multline*}

\textbf{Remark.} One can verify the continuity of the numerator $\psi_1(y)-2\psi_2(y)+2\psi_4(y)-\psi_5(y)$ of the integrand in (\ref{eq:sigmaequal2}) in the neighborhood
of the branch cut $[1,\infty)$ by direct application of Theorem~\ref{th:jump}. Indeed, according to (\ref{eq:jump}) the jump of $\psi_m(y)$ over the ray $y>1$ equals
$$
2{\pi}iy^{m+2}\frac{(\a)_m\Gamma(\Delta(\b+m,3),\Delta(m+3,3))}{(3)_m(\b)_m\Gamma(3,\Delta(\a+m,3))}
G^{3p+4,0}_{3p+4,3p+4}\left(\frac{1}{y^3}\,\vline\begin{array}{c}1,\Delta(\b+m,3),\Delta(m+3,3)\\3,\Delta(\a+m,3)\end{array}\!\!\right).
$$
By the Gauss multiplication formula $\Gamma(\Delta(x+m,3))/(x)_m=2\pi\Gamma(x)3^{1/2-x-m}$, so that the expression
$$
\frac{(\a)_m\Gamma(\Delta(\b+m,3),\Delta(m+3,3))}{(3)_m(\b)_m\Gamma(3,\Delta(\a+m,3))}
$$
is independent of $m$.  Hence, the jump of $\psi_1(y)-2\psi_2(y)+2\psi_4(y)-\psi_5(y)$ is equal to zero if
\begin{multline*}
G^{3p+4,0}_{3p+4,3p+4}\left(\frac{1}{y^3}\,\vline\begin{array}{c}1,\Delta(\b+1,3),\Delta(4,3)\\3,\Delta(\a+1,3)\end{array}\!\!\right)
-2yG^{3p+4,0}_{3p+4,3p+4}\left(\frac{1}{y^3}\,\vline\begin{array}{c}1,\Delta(\b+2,3),\Delta(5,3)\\3,\Delta(\a+2,3)\end{array}\!\!\right)
\\
+2y^3G^{3p+4,0}_{3p+4,3p+4}\left(\frac{1}{y^3}\,\vline\begin{array}{c}1,\Delta(\b+4,3),\Delta(7,3)\\3,\Delta(\a+4,3)\end{array}\!\!\right)
-y^4G^{3p+4,0}_{3p+4,3p+4}\left(\frac{1}{y^3}\,\vline\begin{array}{c}1,\Delta(\b+5,3),\Delta(8,3)\\3,\Delta(\a+5,3)\end{array}\!\!\right)=0.
\end{multline*}
Writing $y^{-3}=t$ and using the shifting property of $G$ function \cite[8.2.2.15]{PBM3} the last formula takes the form
\begin{multline*}
G^{3p+4,0}_{3p+4,3p+4}\left(t\,\vline\begin{array}{c}1,\Delta(\b+1,3),\Delta(4,3)\\3,\Delta(\a+1,3)\end{array}\!\!\right)
-2G^{3p+4,0}_{3p+4,3p+4}\left(t\,\vline\begin{array}{c}2/3,\Delta(\b+1,3),\Delta(4,3)\\8/3,\Delta(\a+1,3)\end{array}\!\!\right)
\\
+2G^{3p+4,0}_{3p+4,3p+4}\left(t\,\vline\begin{array}{c}0,\Delta(\b+1,3),\Delta(4,3)\\2,\Delta(\a+1,3)\end{array}\!\!\right)
-G^{3p+4,0}_{3p+4,3p+4}\left(t\,\vline\begin{array}{c}-1/3,\Delta(\b+1,3),\Delta(4,3)\\5/3,\Delta(\a+1,3)\end{array}\!\!\right)=0.
\end{multline*}
To verify this identity write the definition of $G$ function \cite[8.2.1.1]{PBM3} and collect the integrands under single integral sign to get
$$
\frac{\Gamma(s+3)}{\Gamma(s+1)}-2\frac{\Gamma(s+8/3)}{\Gamma(s+2/3)}+2\frac{\Gamma(s+2)}{\Gamma(s)}-\frac{\Gamma(s+5/3)}{\Gamma(s-1/3)}=0
$$
which is true by the shifting property $\Gamma(z+1)=z\Gamma(z)$.

\medskip

The following identity reduces to Newton-Leibnitz formula under substitution $u=1/(x+t)$:
\begin{equation}\label{eq:foffprime}
\phi(1/x)=\phi_0+\int\limits_{0}^{\infty}\frac{\phi'\left(1/(x+t)\right)}{(x+t)^2}dt,
\end{equation}
where $\phi'(z)=d\phi(z)/dz$. It is valid for differentiable $\phi$ under convergence of the integral in (\ref{eq:foffprime}).
Applying (\ref{eq:foffprime}) to the generalized hypergeometric function we obtain for $p\le{q}$:
$$
{_{p}F_{q}}\left(\left.\!\!\begin{array}{c}\a\\\b\end{array}\right|1/x\!\right)
=1+\frac{(\a)}{(\b)}\int\limits_{0}^{\infty}(x+t)^{-2}{_{p}F_{q}}\left(\left.\!\!\begin{array}{c}\a+1\\\b+1\end{array}\right|\frac{1}{x+t}\!\right)dt,
$$
where we have used the standard formula for the derivative of ${_{p}F_{q}}$ \cite[(16.3.1)]{NIST}. This simple relation is not contained in \cite{PBM3}, but \textit{Mathematica} is able to evaluate the integral on the right hand side.

A straightforward consequence of a result due to Thale \cite{Thale} observed by us in \cite[Theorem~13]{KPJMAA} is that the
generalized hypergeometric function ${}_{p+1}F_{p}(z)$ is univalent in $\Re(z)<1$ as long as it belongs to the Stieltjes class.
Univalent functions in the unit disk satisfy a number of distortion theorems outlined, for instance, in \cite[Theorem~7.1]{Dub}.
Below we illustrate the application of these theorems to the generalized hypergeometric function.  We will need the following inequalities
satisfied by the functions $f$ holomorphic and univalent in the unit disk and normalized by $f(0)=f'(0)-1=0$:
\begin{equation}\label{eq:odnolist1}
\frac{|w|}{(1+|w|)^2}\leq |f(w)|\leq \frac{|w|}{(1-|w|)^2},
\end{equation}
\begin{equation}\label{eq:odnolist2}
\frac{1-|w|}{(1+|w|)^3}\leq {|f'(w)|}\leq \frac{1+|w|}{(1-|w|)^3},
\end{equation}
\begin{equation}\label{eq:odnolist3}
\frac{1-|w|}{|w|(1+|w|)}\leq \frac{|f'(w)|}{|f(w)|}\leq
\frac{1+|w|}{|w|(1-|w|)}.
\end{equation}

\begin{theorem}\label{th:estimate}
Suppose  $\a=(a_1,\ldots,a_{p+1})$ and $\b=(b_1,\ldots,b_p)$ are positive vectors such that, $0<a_1\leq{1}$ and $v_{\a_{[1]},\b}(t)\ge0$ on $[0,1]$, where $v_{\a_{[1]},\b}$ is defined in $(\ref{eq:v-defined})$. Then the following inequalities hold in the half plane $\Re{z}<1$\emph{:}
\begin{equation}\label{eq:ineq1}
\frac{2(\a)|z-2||z|}{(\b)(|z-2|+|z|)^2}\leq|_{p+1}F_p(\a;\b;z)-1|\leq \frac{2(\a)|z-2||z|}{(\b)(|z-2|-|z|)^2},
\end{equation}
\begin{equation}\label{eq:ineq2}
\frac{4(|z-2|-|z|)}{(|z-2|+|z|)^3}\leq
|_{p+1}F_p(\a+1;\b+1;z)|\leq \frac{4(|z-2|+|z|)}{((|z-2|-|z|)^3}
\end{equation}
and
\begin{multline}\label{eq:ineq3}
\frac{2(\b)(|z-2|-|z|)}{(\a)|z||z-2|(|z-2|+|z|)}\leq
\left|\frac{_{p+1}F_p(\a+1;\b+1;z)-1}{_{p+1}F_p(\a;\b;z)}\right|\leq
\frac{2(\b)(|z-2|+|z|)}{(\a)|z||z-2|(|z-2|-|z|)}.
\end{multline}
\end{theorem}
\textbf{Proof}. According to \cite[Theorem~13]{KPJMAA} the function $z\mapsto{}_{p+1}F_p\left(\a;\b; z\right)$ is univalent in the half-plane
 $\Re{z}<1$ under the hypotheses of the theorem.  The M\"{o}bius map $z=z(w)=2w/(w-1)$ effects a biholomorphic bijection between the unit disk $|w|<1$ and
 the half-plane $\Re{z}<1$.  Hence, inequalities (\ref{eq:odnolist1})--(\ref{eq:odnolist3}) hold for the function
$$
f(w)=\frac{(\b)}{2(\a)}\left(1-{}_{p+1}F_p\!\!\left(\!\!\begin{array}{l}\a\\\b\end{array}\vline\,\,z(w)\!\right)\right)
$$
satisfying $f(0)=f'(0)-1=0$.  In view of the derivative formula for the generalized hypergeometric function
and substituting $w=z/(z-2)$ we obtain (\ref{eq:ineq1})--(\ref{eq:ineq3}).$\hfill\square$

\textbf{Remark.} It is a classical fact that equality is attained in (\ref{eq:odnolist1})--(\ref{eq:odnolist3}) for the rotations of the Koebe function $K(w)=w/(1+w)^2$, see \cite[Theorem~7.1]{Dub}. The hypergeometric expression for $K(w(z))$ is given by
$$
K(w)=K\left(\frac{z}{z-2}\right)=\frac{z(z-2)}{4(1-z)^2}=\frac{1}{4}\left(1-\frac{1}{(1-z)^2}\right)
=\frac{1}{4}\left(1-{_{2}F_1}\!\left(\!\!\begin{array}{l}1,2\\1\end{array}\vline\,\,z\!\right)\right).
$$
Comparing this expression with $f(w)$ defined in the proof of Theorem~\ref{th:estimate}, we see that the function $g(z)={_{2}F_1}\!\left(1,2;1;z\right)$ is extremal in Theorem~\ref{th:estimate}, so that equality is attained in (\ref{eq:ineq1})--(\ref{eq:ineq3}) for this function and some values of $z$.  On the other hand, $g(z)$ clearly violates the conditions of Theorem~\ref{th:estimate} as $v_{1,2}(t)=t(t-1)<0$ for $t\in(0,1)$.  Since $g(z)$ is univalent, it is clear that conditions of Theorem~\ref{th:estimate} are far from necessary for univalence of ${}_{p+1}F_p(\a;\b;z)$ in $\Re{z}<1$ and for the validity of (\ref{eq:ineq1})--(\ref{eq:ineq3}). Hence, there is a room for relaxing these conditions.

\section{A curious integral evaluation}

In the final section of the paper we present a curious representation for a finite sum of product ratios of gamma functions as an integral of Meijer's $G$ function.  If the dimension $p=1$, the integrand is expressed via the Gauss hypergeometric function as stated in Corollary~\ref{cr:gammasum-int}.
\begin{theorem}\label{th:gammasum-int}
Suppose $\a,\b\in\C^p$ are such that $0<\Re(a_j)<\Re\left(\sum_{j=1}^{p}b_j\right)$ for $j=1,\ldots,p$, $\alpha,\beta>0$ and $m\in\N_0$. Then the following identity holds
\begin{multline}\label{eq:sum-integral}
\sum\limits_{k=0}^{m}\left\{\frac{\Gamma(\a+k)\Gamma(\a+\alpha+\beta+m-k)}{\Gamma(\b+k)\Gamma(\b+\alpha+\beta+m-k)}
-\frac{\Gamma(\a+\alpha+k)\Gamma(\a+\beta+m-k)}{\Gamma(\b+\alpha+k)\Gamma(\b+\beta+m-k)}\right\}
\\
=\int\limits_{0}^{1}G^{p,p}_{2p,2p}\left(x\left|\begin{array}{l}\!\!1-\a-m,\b+\alpha+\beta\!\!\\\!\!\a+\alpha+\beta,1-\b-m\!\!\end{array}\right.\right)
\!\frac{(1-x^{m+1})(1-x^{\alpha})(1-x^{\beta})}{x^{\alpha+\beta+1}(1-x)}dx.
\end{multline}
\end{theorem}
\textbf{Proof.} Suppose $m+n=p$.  We will write $\b=(\mathbf{bt},\mathbf{bb})$, where $\mathbf{bt}$ stands for the first $m$ components of $\b$ (located in the numerator of the integrand, hence $\mathbf{bt}$ for $\b$-top) and $\mathbf{bb}$ for the lasts $q-m$ components (located in the denominator of the integrand, hence $\mathbf{bb}$ for $\b$-bottom).  Similarly, $\a=(\mathbf{at},\mathbf{ab})$ with $n$ components in $\mathbf{at}$ and $p-n$ in $\mathbf{ab}$. Then, according to \cite[8.2.2.7]{PBM3} and \cite{Wolfram}, Meijer's $G$ function $G^{m,n}_{p,p}(z)$ is piecewise analytic with discontinuity on the unit circle $|z|=1$.  In this case, combining \cite[Theorems~12.5.1]{BealsWong} with \cite[Theorems~12.5.2]{BealsWong} or \cite[8.2.2.3]{PBM3} with \cite[8.2.2.4]{PBM3} and  we can write for $x>0$:
\begin{multline}\label{eq:Gppexpansion}
G^{m,n}_{p,p}\left(x\left|\begin{array}{l}\!\!\a\!\!\\\!\!\b\!\!\end{array}\right.\right)
=H(1-x)\sum\limits_{j=1}^{m}A_jx^{b_j}{}_{p}F_{p-1}\!\!\left(\!\!\begin{array}{l}1-\a+b_j\\1-\b_{[j]}+b_j\end{array}\!\vline\,x\!\right)
\\
+H(x-1)\sum\limits_{j=1}^{n}B_jx^{a_j-1}{}_{p}F_{p-1}\!\!\left(\!\!\begin{array}{l}1+\b-a_j\\1+\a_{[j]}-a_j\end{array}\!\vline\,\frac{1}{x}\!\right),
\end{multline}
where
$$
A_j=\frac{\Gamma(\mathbf{bt}_{[j]}-b_j)\Gamma(1-\mathbf{at}+b_j)}{\Gamma(1-\mathbf{bb}+b_j)\Gamma(\mathbf{ab}-b_j)},
~~~B_j=\frac{\Gamma(a_j-\mathbf{at}_{[j]})\Gamma(1+\mathbf{bt}-a_j)}{\Gamma(1+\mathbf{ab}-a_j)\Gamma(a_j-\mathbf{bb})}
$$
and $H(t)=\partial_t\max\{t,0\}$ is the Heaviside function.  Next, according to \cite[2.24.2.1]{PBM3} (see also \cite[Th.~2.2 and 3.3]{KilSaig}):
\begin{multline*}
\frac{\Gamma(\a+k)\Gamma(\a+\alpha+\beta+m-k)}{\Gamma(\b+k)\Gamma(\b+\alpha+\beta+m-k)}=\int\limits_{0}^{\infty}x^{k-1}
G^{p,p}_{2p,2p}\left(x\left|\begin{array}{l}\!\!1-\a-\alpha-\beta-m,\b\!\!\\\!\!\a,1-\b-\alpha-\beta-m\!\!\end{array}\right.\right)dx
\\
=\int\limits_{0}^{\infty}x^{k-\alpha-\beta-1}
G^{p,p}_{2p,2p}\left(x\left|\begin{array}{l}\!\!1-\a-m,\b+\alpha+\beta\!\!\\\!\!\a+\alpha+\beta,1-\b-m\!\!\end{array}\right.\right)dx,
\end{multline*}
where in the last equality we have used the shifting property \cite[8.2.2.15]{PBM3}.  Similarly,
\begin{multline*}
\frac{\Gamma(\a+\alpha+k)\Gamma(\a+\beta+m-k)}{\Gamma(\b+\alpha+k)\Gamma(\b+\beta+m-k)}=\int\limits_{0}^{\infty}x^{k-1}
G^{p,p}_{2p,2p}\left(x\left|\begin{array}{l}\!\!1-\a-\beta-m,\b+\alpha\!\!\\\!\!\a+\alpha,1-\b-\beta-m\!\!\end{array}\right.\right)dx
\\
=\int\limits_{0}^{\infty}x^{k-\beta-1}
G^{p,p}_{2p,2p}\left(x\left|\begin{array}{l}\!\!1-\a-m,\b+\alpha+\beta\!\!\\\!\!\a+\alpha+\beta,1-\b-m\!\!\end{array}\right.\right)dx.
\end{multline*}
Note that the $G$ functions in the integrand are the same, so that we can also write the left hand side of (\ref{eq:sum-integral}) as
\begin{equation}\label{eq:gamasumtoGfunction}
\text{LHS of (\ref{eq:sum-integral})}=\int\limits_{0}^{\infty}x^{k-\beta-1}(x^{-\alpha}-1)
G^{p,p}_{2p,2p}\left(x\left|\begin{array}{l}\!\!1-\a-m,\b+\alpha+\beta\!\!\\\!\!\a+\alpha+\beta,1-\b-m\!\!\end{array}\right.\right)dx.
\end{equation}
However, the function in the integrand is piecewise analytic, so essentially we are dealing with a sum of two different integrals. To reduce it to a single integral we
apply (\ref{eq:Gppexpansion}) to the $G$ function in the integrand:
\begin{multline}\label{eq:Gpp2p2p}
G^{p,p}_{2p,2p}\left(x\left|\begin{array}{l}\!\!1-\a-m,\b+\alpha+\beta\!\!\\\!\!\a+\alpha+\beta,1-\b-m\!\!\end{array}\right.\right)
\\
=H(1-x)\sum\limits_{j=1}^{p}M_jx^{a_j+\alpha+\beta}{}_{2p}F_{2p-1}\!\!\left(\!\!\begin{array}{l}\a+m+\alpha+\beta+a_j,1-\b+a_j\\1-\a_{[j]}+a_j,\b+m+\alpha+\beta+a_j\end{array}\!\vline\,x\!\right)
\\
+H(x-1)\sum\limits_{j=1}^{p}M_jx^{-a_j-m}{}_{2p}F_{2p-1}\!\!\left(\!\!\begin{array}{l}\a+m+\alpha+\beta+a_j,1-\b+a_j\\1-\a_{[j]}+a_j,\b+m+\alpha+\beta+a_j\end{array}\!\vline\,\frac{1}{x}\!\right),
\end{multline}
where
\begin{equation}\label{eq:Mj}
M_j=\frac{\Gamma(\a_{[j]}-a_j)\Gamma(\a+m+\alpha+\beta+a_j)}{\Gamma(\b+m+\alpha+\beta+a_j)\Gamma(\b-a_j)}.
\end{equation}
Note the surprising fact that both the constants and the hypergeometric functions in the two sums are exactly the same. Next, substitute this expansion into the integrand in (\ref{eq:gamasumtoGfunction}), carry out change of variable $x\to1/x$ in the second integral and add up the resulting integrals to get
\begin{multline}\label{eq:sum-to-int-of-pFp-1}
\text{LHS of (\ref{eq:sum-integral})}=\int\limits_{0}^{1}\!\biggl[\!\sum_{j=1}^{p}M_jx^{a_j-1}
\\
\times{}_{2p}F_{2p-1}\!\!\left(\!\!\begin{array}{l}\a+a_j+\alpha+\beta+m,1-\b+a_j\\1-\a_{[j]}+a_j,\b+a_j+\alpha+\beta+m\end{array}\!\vline\,x\!\right)\biggr]
\!\frac{(1-x^{m+1})(1-x^{\alpha})(1-x^{\beta})}{1-x}dx.
\end{multline}
In order to obtain (\ref{eq:sum-integral}) it remains to observe that the expression in brackets equals the left hand side of (\ref{eq:Gpp2p2p}) times $H(1-x)x^{-\alpha-\beta}$.
Conditions for convergence follow by applying the expansion \cite[(2.19)]{KPSigma} of ${}_{2p}F_{2p-1}(z)$ in the neighborhood of unity (due to N{\o}rlund and B\"{u}hring) to the first sum in (\ref{eq:Gpp2p2p}).  $\hfill\square$

Taking $p=1$ in the above theorem we immediately get the next
\begin{corollary}\label{cr:gammasum-int}
Suppose $\Re(b+1)>\Re(a)>0$, $\alpha,\beta>0$ and $m\in\N_0$. Then the following identity holds
\begin{multline}\label{eq:sum-integral1}
\sum\limits_{k=0}^{m}\left\{\frac{\Gamma(a+k)\Gamma(a+\alpha+\beta+m-k)}{\Gamma(b+k)\Gamma(b+\alpha+\beta+m-k)}
-\frac{\Gamma(a+\alpha+k)\Gamma(a+\beta+m-k)}{\Gamma(b+\alpha+k)\Gamma(b+\beta+m-k)}\right\}
\\
=\frac{\Gamma(2a+\alpha+\beta+m)}{\Gamma(a+b+\alpha+\beta+m)\Gamma(b-a)}
\\
\times\!\!\int\limits_{0}^{1}\!\!x^{a-1}(1-x^{m+1})(1-x^{\alpha})(1-x^{\beta})
{_2F_1}\!\!\left(\!\!\begin{array}{l}2a+\alpha+\beta+m,1+a-b\\a+b+\alpha+\beta+m\end{array}\vline\,\,x\!\right)\frac{dx}{1-x}.
\end{multline}
\end{corollary}

\paragraph{Acknowledgements.} This research has been supported by the Russian Science Foundation under project 14-11-00022.

\end{document}